\documentclass[11pt]{amsart}

\usepackage{pdfsync}

\usepackage{amssymb,graphics,epsfig,overpic,amscd,amsmath}
\usepackage[usenames,dvipsnames]{color}
\usepackage[colorlinks=true,linkcolor=blue,citecolor=BrickRed]{hyperref}

\newtheorem{thm}{Theorem}[section]

\newtheorem{condition}[thm]{Condition}
\newtheorem{question}[thm]{Question}

\theoremstyle{definition}

\newcommand{\R}{\mathbf{R}}

\newcommand{\T}{\mathbf{T}}
\newcommand{\ol}{\overline}
\newcommand{\C}{\mathcal{C}}

\renewcommand{\S}{\mathbf{S}}

\renewcommand{\tilde}{\widetilde}
\renewcommand{\epsilon}{\varepsilon}

\DeclareMathOperator{\Cone}{Cone}

\title[]{Directed immersions of closed manifolds}

\author{Mohammad Ghomi}
\address{School of Mathematics, Georgia Institute of Technology,
Atlanta, GA 30332}
\email{ghomi@math.gatech.edu}
\urladdr{www.math.gatech.edu/$\sim$ghomi}

\date{\today \,(Last Typeset)}
\thanks{The research of the author was supported in part by NSF grant DMS-0806305.}

\begin{document}

\begin{abstract}
Given any finite subset  $X$ of the sphere $\S^{n}$, $n\geq 2$, which includes no pairs of antipodal points, we  explicitly construct smoothly immersed closed orientable hypersurfaces in Euclidean space $\R^{n+1}$  whose Gauss map misses $X$. In particular, this answers a question of M. Gromov.
\end{abstract}

\maketitle

\section{Introduction}
To every $(\C^1)$ immersion $f\colon M^n\to\R^{n+1}$ of a closed oriented $n$-manifold $M$, there corresponds  a unit normal vector field or \emph{Gauss map} $G_f\colon M\to\S^n$, which generates a set  $G_f(M)\subset\S^n$ known as  the \emph{spherical image} of $f$. Conversely, one may ask, c.f.
\cite[p. 3]{gromov:questions}:  \emph{for which sets $A\subset\S^n$ is there an immersion $f\colon M\to \R^{n+1}$  such that $G_f(M)\subset A$?} Such a mapping would be called an \emph{$A$-directed immersion} of $M$ \cite{eliashberg&mishachev,gromov:pdr,rs:compression,spring:directed}. It is well-known that when $A\neq\S^n$,  $f$  must have double points (Note \ref{note:embedding}), and $M$ must be parallelizable, e.g., $M$ can only be the torus $\T^2$ when $n=2$ (Note \ref{note:torus}). Furthermore, the only known necessary condition on $A$ is the elementary observation that $A\cup -A=\S^2$, while there is also a sufficient condition due to  Gromov \cite[Thm. $(D')$,  p. 186]{gromov:pdr}:

\begin{condition}\label{cond}
$A\subset\S^n$ is open, and there is a point $p\in\S^n$ such that the intersection of $A$ with  each great circle passing through $p$ includes a (closed) semicircle.
 \end{condition}
 
 A \emph{great circle}  is the intersection of $\S^n$ with a $2$-dimensional subspace   of $\R^{n+1}$.
 Note  that, when $n\geq 2$, examples of sets $A\subset\S^n$ satisfying the above condition include those which are the complement of a finite set of points without  antipodal pairs. Thus the spherical image of a closed hypersurface can be remarkably flexible. Like most $h$-principle or convex integration type arguments, however, the proof  does not yield specific examples. It is therefore natural to ask, for instance:
 
 \begin{question}[\cite{gromov:pdr}, p. 186] \label{question}
``Is there a `simple' immersion $\T^2\to\R^3$ whose spherical image misses the four vertices of a regular tetrahedron in $\S^2$?" 
\end{question}

\begin{figure}[h] 
   \centering
   \includegraphics[height=1.1in]{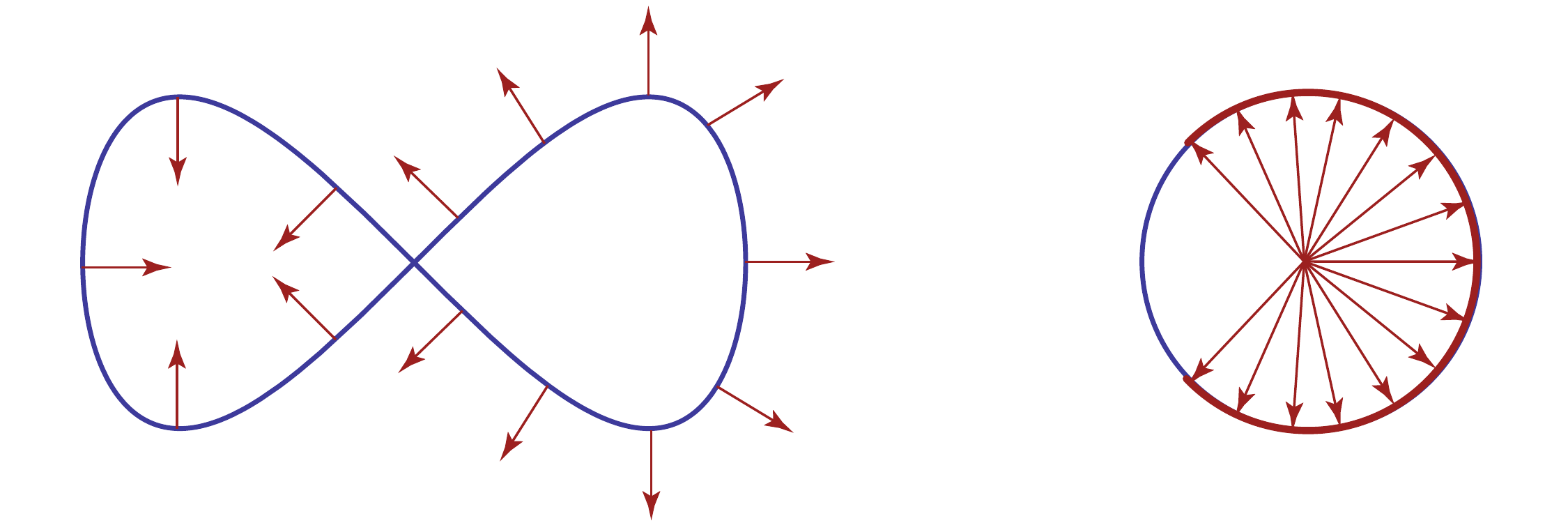} 
   \caption{}
   \label{fig:82}
\end{figure}

Here we give an affirmative answer to this question (Section \ref{sec:example}), and more generally present a short  constructive proof of the sufficiency of  a slightly stronger version of Condition \ref{cond} for the existence of $A$-directed immersions of parallelizable manifolds $M^{n-1}\times\S^1$,
where $M^{n-1}$ is  closed and orientable.  Any such manifold  admits an immersion $f\colon M^{n-1}\to\R^n\times\{0\}\subset\R^{n+1}$ (Note \ref{note:immersion}). We then extend $f$ to $M^{n-1}\times\S^1$ by using the
 \emph{figure eight curve} 
\begin{equation}\label{eq:eight}
E_\delta(t):=\big(\cos(t),\,\delta\sin(2t)\big)
\end{equation}
to put a copy of $\S^1\simeq\R/2\pi$ in each normal plane of $f$, as described below. Note that the midpoint of  $G_{E_\delta}(\S^1)$ is assumed to be at $(1,0)$; see Figure \ref{fig:82} which shows $E_{1/2}$ and its spherical image. Further, the unit normal bundle of $f$ may be naturally identified with the pencil of great circles of $\S^n$ passing through $(0,\dots, 0,1)$.

 \begin{thm}\label{thm:main}
 Let $A\subset\S^n$ satisfy Condition \ref{cond} with respect to $p=(0,\dots,0,1)$. Further, if $n\geq 3$, suppose that the semicircle in Condition \ref{cond} contains $p$, or that no great circle through $p$ is contained in $A$. Let
 $f\colon M^{n-1}\to\R^{n}\times\{0\}\subset\R^{n+1}$ be a smooth ($\C^\infty$) immersion of a closed orientable $(n-1)$-manifold,
 and, for every $q\in M$, $C_q\subset \S^n$ be the unit normal space of $f$ at $q$. Then there is a
 smooth orthogonal frame $\{N_i\colon M\to\S^n\}$, $i=1$, $2$, for the normal bundle of $f$
 such that the semicircle in $C_q$ centered at $N_1(q)$ lies in $A$.    For any such frame, and sufficiently small  $\epsilon$, $\delta>0$,
\begin{equation}\label{eq:f}
 F(q,t):=f(q)+\epsilon\sum_{i=1}^2E^i_\delta(t)N_i(q)
 \end{equation}
 yields a smooth $A$-directed immersion $M\times\S^1\to\R^{n+1}$, where $E_\delta^i$ are the components of the figure eight curve  $E_\delta$ given by \eqref{eq:eight}. 
  \end{thm}
 
It is not known if Condition \ref{cond} is necessary for the existence of $A$-directed  closed hypersurfaces, and the question posed in the first paragraph is open, even for $n=2$. 
See \cite{ghomi:shadow,ghomi:tangents, ghomi&tabachnikov} for some other recent results on Gauss maps of closed submanifolds, 
\cite{ghomi:topology,hartman&nirenberg, milnor:immersion, wu:spherical} for still more studies of spherical images, and \cite{spring:history} for historical  background.
 %%%%%%%%%%%%%%%%%%%%%%%%%%%%%%%%%%%%%%%%%%%%
\section{Example}\label{sec:example}
%%%%%%%%%%%%%%%%%%%%%%%%%%%%%%%%%%%%%%%%%%%%

\begin{figure}[h] 
    \centering
    \includegraphics[height=1.4in]{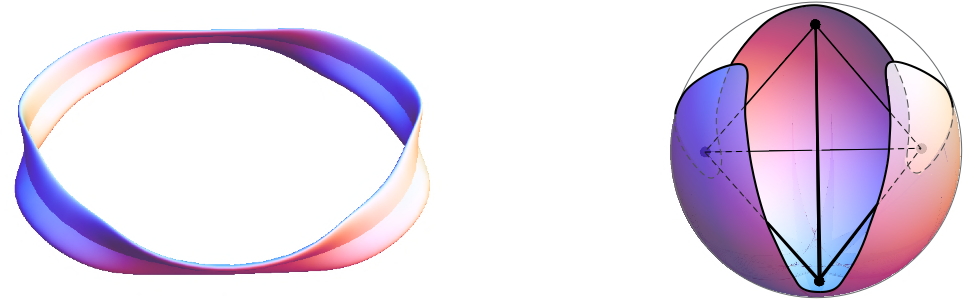} 
    \caption{}
    \label{fig:ribbonandsphere}
 \end{figure}

If $A=\S^2\smallsetminus X$ for a finite set $X$ without antipodal pairs,  we may always find a point $p\in\S^2$ with respect to which $A$ satisfies the hypothesis of Theorem \ref{thm:main}  (e.g., let $p\not\in X$ be  in the complement of all great circles which pass through at least two points of $X$ other than $-p$).  After a rigid motion (which may be arbitrarily small) we may  assume that 
 $p=(0,0,1)$ or $(0,0,-1)$, and let $f(\theta):=(\cos(\theta),\sin(\theta),0)$ be the standard immersion of $\S^1\simeq\R/2\pi$ in $\R^3$. Then the desired framing for the normal bundle of $f$ may always take the form
\begin{equation}\label{eq:z}
N_1(\theta):=f'(\theta)\times N_2(\theta), \quad\quad\quad N_2(\theta):=\frac{\big(\cos(\theta),\,\sin(\theta),\,z(\theta)\big)}{\sqrt{1+z^2(\theta)}},
\end{equation}
where $z\colon\R/2\pi\to \R$ is a smooth function with $z(\theta)=-z(\theta+\pi)$  and  such that $X$ is contained entirely in one of the components of $\S^2-N_2(\S^1)$. For instance, when $X$ is the vertices of a regular tetrahedron, we may set 
$
z(\theta):=\cos(3\theta)
$
in \eqref{eq:z}. Then,  for $\epsilon$, $\delta\leq 1/8$, the mapping $F(\theta,t)$ given by \eqref{eq:f}  yields an  immersion $\T^2\simeq \R/2\pi\times\R/2\pi\to\R^3$ which answers Question \ref{question}.
The resulting surface, for $\epsilon=\delta=1/8$, is depicted in  Figure \ref{fig:ribbonandsphere}   
together with its spherical image (note that here $p=(0,0,-1)$).  
To find $z(\theta)$ in  general, we may order the points in $X'\cup -X'$, where $X':=X\smallsetminus\{-p\}$, according to  their ``longitude"  $\theta$, and connect them by geodesic segments to obtain a  simple closed symmetric curve $\gamma(\theta)$. A perturbation of $\gamma$ then yields a smooth symmetric curve $\tilde\gamma$ such that $X$ is contained in one of the components of $\S^2-\tilde\gamma(\S^1)$. The third coordinate of $\tilde\gamma$ gives our desired height function $z$.
  
 %%%%%%%%%%%%%%%%%%%%%%%%%%%%%%%%%%%%%%%%%%%%%%%% 
 \section{Proof of Theorem \ref{thm:main}}\label{sec:proof}
 %%%%%%%%%%%%%%%%%%%%%%%%%%%%%%%%%%%%%%%%%%%%%%%%
\subsection{}\label{subsec:proof1}
First we construct the frame  $\{N_i\}$. 
For every $q\in M$,  $C_q$ is a great circle passing through $p$. So it contains a  semicircle in $A$ by assumption (Condition \ref{cond}). Let $m_{q}\subset C_q$ be the set of midpoints of all such  semicircles. We need  to find a smooth map $N_1\colon M\to\S^n$ such that  $N_1(q)\in m_{q}$ for all $q\in M$.
To this end note that $m_{q}$ is  open and connected. Further, if $m_{q}$ contains any pairs of antipodal points, then $m_q=C_q$; otherwise,  $m_{q}$ lies in the interior a semicircle of $C_{q}$. Consequently,  
$$
\Cone(m_q):=\{\,\lambda x\mid x\in m_q, \;\text{and}\;\;\lambda\geq 0\,\},
$$
 is a convex set in $\R^{n+1}$. In particular, for any finite set of points $x_i\in \Cone(m_q)$ and numbers $\lambda_i\geq 0$, $\sum_i\lambda_i x_i\in \Cone(m_q)$.
Now let $B$ be the set of all points $q\in M$ such that $m_{q}\neq C_q$. Then $B$ is closed (and  therefore compact) since $M\smallsetminus B$ is open; indeed the set of great circles contained in $A$ is open, since $A$ is open. Further note that for any point $q\in M$, normal vector $x\in m_{q}$,  and continuous local extension $v$ of $x$ to a normal vector field of $M$, we have $v(q')\in m_{q'}$ for all $q'$ in an open neighborhood $U$ of $q$  (because the set of semicircles contained in $A$ is open). Let $\{v_i\colon U_i\to\S^n\}$, $i=1,\dots,k$, be a finite collection of such local vector fields  so that $\cup_iU_i$ covers $B$ and $v_i$ are smooth.  Also let $\{\phi_i\colon M\to\R\}$ be a smooth partition of unity subordinate to $\{U_i\}$, and, for $q\in \cup_i U_i$, set
$$
N_1(q):=\frac{\sum_{i=1}^k \phi_i(q) v_i(q)}{\|\sum_{i=1}^k \phi_i(q) v_i(q)\|}.
$$
 If $q\in B$, then $v_i(q)\in m_{q}$ which lies in the interior of a semicircle $S\subset C_q$, and so $\|\sum_{i=1}^k \phi_i(q) v_i(q)\|\neq 0$. Indeed, if $x$ is  the midpoint of $S$, then
 $
 \langle \sum_{i=1}^k \phi_i(q) v_i(q),x\rangle=\sum_{i=1}^k \phi_i(q)\big\langle  v_i(q),x\big\rangle>0.
 $
  Thus  $N_1$ is well defined (and smooth) on an open neighborhood $V$ of $B$. Further, $N_1(q)\in m_{q}$, for all $q\in V$, since $\Cone(m_q)$ is convex. In particular we are done if $B=M$; otherwise, note that  we may write 
  \begin{equation}\label{eq:n1}
  N_1(q)=\cos\big(\theta(q)\big) \,p+\sin\big(\theta(q)\big)\,G_f(q),
  \end{equation}
  for some function $\theta\colon V\to \R$, since $G_f$ is well defined due to the orientability of $M$, and thus $\{p, G_f(q)\}$ forms an orthonormal basis for the normal plane $df(T_qM)^\perp$. 
  Further, it is easy to see that we may choose $\theta$ continuously (and therefore smoothly) if $n=2$. This also holds for $n>2$ if each $C_q$ contains a semicircle passing through $p$; for then  $\theta$ is uniquely determined within  the range $[-\pi/2,\pi/2]$. Indeed, we may choose the vectors $v_i$ above so that $\langle v_i(q), p\rangle\geq 0$ which would in turn yield that $\langle N_1(q), p\rangle \geq 0$. 
  Now  let $V'$ be an open neighborhood of $B$ with  closure $\ol {V'}\subset V$. 
  Using Tietze's theorem, followed by a perturbation and a gluing, we may extend $\theta|_{V'}$ smoothly to all of $M$. Then \eqref{eq:n1}  yields the desired vector field on $M$, since for any $q\in M\smallsetminus B$, $N_1(q)\in C_q=m_q$. Finally, set 
  $$
  N_2(q):= \sin\big(\theta(q)\big)\,p -\cos\big(\theta(q)\big)\,G_f(q).
  $$
  
\subsection{}\label{subsec:proof2}
It remains to show that  $G_F(M\times\S^1)\subset A$, for small $\epsilon$, $\delta>0$.  For all $q\in M$,  $C_q\cap A$ contains an arc of length $\geq\pi+\alpha$ with midpoint $N_1(q)$ for some uniform constant $\alpha > 0$. 
Indeed, if  we let $g(q)$ be the supremum of lengths of  all arcs in $C_q\cap A$ with midpoint $N_1(q)$, then $g\colon M\to\R$ is lower semicontinuous, i.e., $\lim_{q\to q_0}g(q)\geq g(q_0)$, since  $A$ is open. Thus, since  $g>\pi$ and $M$ is compact, $g\geq\pi+\alpha$.
Now  choose $\delta>0$ so small that the length $\ell$ of the spherical image of  $E_\delta$  is $\leq \pi+\alpha$ (this is possible since $\ell\to \pi$ as $\delta\to 0$). Next, for  $(q,t)\in M\times\S^1$,  let  $\tilde G_F(q,t)$ be the normalized projection of $G_F(q,t)$ into  $df(T_q M)^\perp$, i.e., 
$$
\tilde G_F(q,t):=\frac{\sum_{i=1}^2 \Big\langle G_F(q,t),N_i(q)\Big\rangle N_i(q)}{\sqrt{\sum_{i=1}^2\Big\langle G_F(q,t),N_i(q)\Big\rangle^2}}.
$$
Also, for fixed $t\in\S^1$, let  $F_t(q):=F(q,t)$. Then, by the tubular neighborhood theorem, $F_t\colon M\to\R^{n+1}$ is a smooth immersion for small $\epsilon$. Further, as $\epsilon\to 0$, $F_t$ converges to $f$ with respect to the $\C^1$-topology. Thus, for each $q\in M$, the normal plane  $dF_t(T_qM)^\perp$ (which contains $G_F(q,t)$) converges to  $df(T_q M)^\perp$. Consequently $G_F$ is well-defined for small $\epsilon$, and converges to $\tilde G_F$ as $\epsilon\to 0$.
So it suffices to check that $\tilde G_F(M\times\S^1)\subset A$, which follows from our choice of $\delta$. Indeed
 for each  fixed $q\in M$, $\tilde G_F(\{q\}\times\S^1)$ is the spherical image of  the figure eight curve $\sum_{i=1}^2 E^i_\delta(t)N_i(q)$ in $df(T_q M)^\perp$, which is an arc of $C_{q}$  with midpoint $N_1(q)$ and length $\leq\pi+\alpha$. \qed

%%%%%%%%%%%%%%%%%%%%%%%%%%%%%%%%%%%%%%%
 \section{Notes}
 %%%%%%%%%%%%%%%%%%%%%%%%%%%%%%%%%%%%%%%%
 
\subsection{}\label{note:embedding}

It is well-known that 
$G_f(M)=\S^n$ for any embedding $f\colon M^n\to\R^{n+1}$ of a closed oriented $n$-manifold \cite[p. 187]{gromov:pdr}. More generally, this also holds   for 
``Alexandrov embeddings", i.e., immersions $f\colon M\to\R^{n+1}$ which may be extended to an immersion $\ol f\colon \ol M\to\R^{n+1}$ of a compact ($n+1$)-manifold $\ol M$ with $\partial\ol M=M$. Indeed if $v$ is any vector field along $M$ which points ``outward" with respect to $\ol M$, then for $p\in M$, the normalized projection of $df(v(p))$ into the line $df(T_p M)^\perp$ defines a normal vector field $M\to\S^n$ which coincides with $G_f$ (after a reflection of $G_f$ if necessary).  Then, for any  $u\in\S^n$, if  $p$ is a point which maximizes the height function $\langle \cdot,u\rangle$ on $M$, we have $G_f(p)=u$. On the other hand, being only regularly  homotopic to an embedding, is not enough to ensure that $G_f(M)=\S^n$. Indeed the example in Figure \ref{fig:ribbonandsphere} is regularly homotopic to  an embedded torus of revolution \cite{pinkall:immersion}.

\subsection{}\label{note:torus}
If  $G_f(M)\neq \S^n$ for an immersion $f\colon M^n\to\R^{n+1}$ of an oriented $n$-manifold, then, as is well-known \cite{milnor:immersion}, $M$ must be parallelizable. Here we include a brief geometric argument for this fact. If $(0,\dots,0,1)\not\in G_f(M)$, we may define a continuous map  $F\colon TM\to\R^n\simeq\R^n\times\{0\}\subset\R^{n+1}$ as follows, c.f. \cite[Lemma 2.2]{ghomi&kossowski}. There is a continuous map $\S^n\smallsetminus \{(0,\dots,0,1)\}\overset{\rho}{\to}SO(n+1)$, $u\mapsto \rho_u$  such that $\rho_u(u)=(0,\dots,0,-1)$. Let $\pi\colon TM\to M$ be the canonical projection, and for $X\in TM$ set $F(X):=\rho_{G_f(\pi(X))}(df(X))$. Also let $F_p:=F|_{T_pM}$. Then $\{F_p^{-1}(e_i)\}$, where $\{e_i\}$ is a fixed basis of $\R^{n}$, gives a framing for $TM$ as desired.  So in particular, when $M$ is closed and $n=2$, we have $M=\T^2$. The last observation also follows from Gauss-Bonnet theorem via degree theory when $f$ is $\C^2$; since if $G_f(M)\neq\S^2$, then
$$
0=\text{deg}(G_f)=\frac{1}{4\pi}\int_M\det(dG_f)=\frac{1}{4\pi} \int_M K=\frac{1}{2}{\chi(M)},
$$
where $K$  is the Gaussian curvature and $\chi$ is the Euler characteristic. 

\subsection{}\label{note:immersion}
To generate some concrete examples of the immersions  $f\colon M^{n-1}\to\R^n\simeq\R^{n}\times\{0\}$ in Theorem \ref{thm:main}, note that if $f_0\colon M_0^{n-k-1}\to\R^{n-k}\times\{0\}$ is any immersion such that $f_0(M_0)$ is disjoint from the subspace $L:=\R^{n-k-1}\times\{(0,0)\}$, then spinning $f_0$ about $L$ yields an immersion $f_1\colon M_0\times\S^1\to\R^{n-k+1}$ given by
$$
f_1(q,t):=
\left[\begin{array}{ccc|c}
&\text{\large I}&
 & \text{\large $0$}\\ 
 & &  &\\
 \hline
&\text{\large $0$}& &
\begin{array}{rl}
\vspace{-0.5em}
& \\
\cos(t) &  \sin(t)\\
-\sin(t) & \cos(t)\\
\end{array}
\end{array}\right]
\left[\begin{array} {c}
f^1_0(q)\\ 
\vdots\\
f^{n-k}_0(q)\\
 0\\
 \end{array}\right],
$$
 where $f_0^i$ are the components of $f_0$. Thus, for instance, one may inductively construct  immersions of $\S^{n-k-1}\times\T^k$ in $\R^n$, for  $k=1,\dots, n-1$. More generally, if $M^{n-1}\times \S^1$ is parallelizable, then so is the open manifold $M^{n-1}\times(0,1)$, which may be immersed in $\R^n$ \cite{hirsh:immersion} by the $h$-principle \cite{gromov:pdr},  or more specifically, the ``holonomic approximation theorem" of Eliashberg-Mishachev \cite{eliashberg&mishachev, ghomi&kossowski}.

\section*{Acknowledgements}
The author thanks Misha Gromov for his interesting question in \cite[p. 186]{gromov:pdr}, and David Spring who first called the author's attention to that problem and pointed out a correction in an earlier draft of this work.

\bibliographystyle{abbrv}
\bibliography{references}

\end{document}